\documentclass[11pt,reqno]{amsart}

\usepackage{amscd,amssymb,amsmath,amsthm}
\usepackage[arrow,matrix]{xy}
\usepackage{graphicx}
\usepackage{epstopdf}
\usepackage{color}
\newtheorem{thm}{Theorem}
\newtheorem{defn}{Definition}
\newtheorem{lemma}{Lemma}
\newtheorem{pro}{Proposition}
\newtheorem{rk}{Remark}

\numberwithin{equation}{section} 

\def\Q{\mathbb Q}

\def\R{\mathbb{R}}

\def\N{\mathbb N}
\def\Z{\mathbb{Z}}

\begin{document}

\title[Group structure of the $p$-adic ball...]{Group structure of the $p$-adic ball and dynamical system of isometry on a sphere}

\author{I.~A.~Sattarov}

\address{I.~A.~Sattarov}

\address{1. Namangan State University, 316, Uychi
 str. Namangan, 160100, Uzbekistan.}

\address{2. Institute of mathematics,
9, University str. Tashkent, 100174, Uzbekistan,}

\email {sattarovi-a@yandex.ru}

\begin{abstract}
In this paper the group structure of the $p$-adic ball and sphere are studied. The dynamical system of isometry defined on invariant sphere is investigated.
We define the binary operations $\oplus$ and $\odot$ on a ball and sphere respectively, and prove that this sets are compact topological abelian group with respect to the operations. Then we show that any two balls (spheres) with positive radius are isomorphic as groups.
We prove that the Haar measure introduced in $\Z_p$ is also a Haar measure on an arbitrary balls and spheres.
We study the dynamical system generated by the isometry defined on a sphere, and
show that the trajectory of any initial point that is not a fixed point isn't convergent.
We study ergodicity of this $p$-adic dynamical system with respect to normalized Haar measure
reduced on the sphere. For $p\geq 3$ we prove that the dynamical
systems are not ergodic. But for $p=2$ under some conditions the dynamical system may be ergodic.
\end{abstract}

 \maketitle

\section{Introduction}

In 1916, Ostrovsky proved that any non-trivial norm defined in the field of rational numbers is equivalent to either the usual absolute value or the $p$-adic norm for some prime number $p$. This means that it is enough to look at "two types" of norms (ordinary absolute value and $p$-adic norm) in the field of rational numbers. As a result was concluded that when it is impossible to explain natural processes through one of them, the other one comes to help, and they are not better or worse than each other, but complement each other.

Currently, $p$-adic theory, and more broadly ultrametric analysis, is a rapidly developing area that finds applications in various sciences (physics, biology, genetics, cognitive sciences, computer science, informatics, cryptology, numerical methods, etc.). To the state-of-the-art the interested reader is referred to, e.g., \cite{ARS}-\cite{Wal} and references therein.

In mathematics, an invariant measure is a measure that is preserved by some function.
Ergodic theory is the study of invariant measures in dynamical systems.
Also, this theory often deals with ergodic transformations.
The intuition behind such transformations, which act on a given set, is that they
do a thorough job of "mixing" the elements of this set. Here is the formal definition.

\begin{defn}\label{erg}\cite{Wal}
Let $T:X\to X$ be a measure-preserving transformation on a measure space $(X, \Sigma, \mu)$, with $\mu(X) = 1$.
Then $T$ is ergodic if for every $E$ in $\Sigma$ with $T^{-1}(E) = E$, either $\mu(E) = 0$ or $\mu(E) = 1$.
\end{defn}

In \cite{An1}-\cite{An4}, \cite{M1}-\cite{M3}, on ergodicity of $p$-adic dynamical systems are considered in the set of $p$-adic integers $\Z_p$.
One of the main reasons for this is the Haar measure introduced in $\Z_p$ as an analogue of the Lebesgue measure on the segment $[0,1]\subset\R$, and ergodicity is studied with respect to this Haar measure.

The questions about the ergodicity of the dynamics considered on invariant spheres with different centers and different radii remain open, because the Haar measure is defined in $\sigma$-algebra composed of subsets of a set with an algebraic structure, and is invariant with respect to the operation. The following questions are interesting: what algebraic structure exists in the considered spheres?; is it possible to determine the Haar measure in these spheres? In this work we will answer such questions.

The paper is organized as follows: in
Section 2 we give some preliminaries. In
Section 3 we define the operations $\oplus$ and $\odot$ on a ball and sphere in $\Q_p$, respectively, and prove that this sets are compact topological abelian group. Then we show that any two balls (spheres) with positive radius are isomorphic as groups.
In Section 4 we prove that the Haar measure introduced in $\Z_p$ is also a Haar measure on an arbitrary balls and spheres.
In Section 5 we study the dynamical system generated by the isometry defined on the invariant sphere, and
we show that the trajectory of any initial point that is not a fixed point isn't convergent.
We study ergodicity of this $p$-adic dynamical system with respect to normalized Haar measure
reduced on the sphere. For $p\geq 3$ we prove that the dynamical
systems are not ergodic. But for $p=2$ under some conditions the dynamical system may be ergodic.

\section{Preliminaries}

Let $\Q$ be the field of rational numbers. The greatest common
divisor of the positive integers $n$ and $m$ is denotes by
$(n,m)$. Every rational number $x\neq 0$ can be represented in the
form $x=p^r\frac{n}{m}$, where $r,n\in\mathbb{Z}$, $m$ is a
positive integer, $(p,n)=1$, $(p,m)=1$ and $p$ is a fixed prime
number.

The $p$-adic norm of $x$ is given by
$$
|x|_p=\left\{
\begin{array}{ll}
p^{-r}, & \ \textrm{ for $x\neq 0$},\\[2mm]
0, &\ \textrm{ for $x=0$}.\\
\end{array}
\right.
$$
It has the following properties:

1) $|x|_p\geq 0$ and $|x|_p=0$ if and only if $x=0$,

2) $|xy|_p=|x|_p|y|_p$,

3) the strong triangle inequality
$$
|x+y|_p\leq\max\{|x|_p,|y|_p\},
$$

3.1) if $|x|_p\neq |y|_p$ then $|x+y|_p=\max\{|x|_p,|y|_p\}$,

3.2) if $|x|_p=|y|_p$ then $|x+y|_p\leq |x|_p$,

this is a non-Archimedean one.

The completion of $\Q$ with  respect to $p$-adic norm defines the
$p$-adic field which is denoted by $\Q_p$ (see \cite{Ko}).

For any $a\in\Q_p$ and $r>0$ denote
$$
U_r(a)=\{x\in\Q_p : |x-a|_p<r\},\ \ V_r(a)=\{x\in\Q_p :
|x-a|_p\leq r\},
$$
$$
S_r(a)=\{x\in\Q_p : |x-a|_p= r\}.
$$

\section{Group structure of the $p$-adic ball and sphere}

Let $\Q_p$ be a field of $p$-adic numbers and $a\in\Q_p$. We consider the ball $V_r(a)$ and the sphere $S_r(a)$ with radius $r\in\{p^m: \, m\in \Z\}$ and the center at $a$.

We define the binary operations $\oplus$ and $\odot$ by the following way
\begin{equation}\label{oplus}
x\oplus y=x+y-a, \ \ \ \ \ \ \ \ \ \ \ \ \ \ \ \ \forall x,y\in V_r(a),
\end{equation}
\begin{equation}\label{odot}
x\odot y=r(x-a)(y-a)+a, \ \ \ \ \forall x,y\in S_r(a).
\end{equation}

\begin{thm}\label{gr}
\begin{itemize}
\item[a)] The pair $(V_r(a),\oplus)$ is a compact topologically additive abelian group.
\item[b)] The pair $(S_r(a),\odot)$ is a compact topologically multiplicative abelian group.
\end{itemize}
\end{thm}

\begin{proof}
It is known \cite{Kat} that any ball and sphere are compact sets in $\mathbb Q_p$. It can be easily seen from (\ref{oplus}) and (\ref{odot}) that the operations $\oplus$ and $\odot$ are commutative.

Let $x,y,z\in V_r(a)$. Then we have $(x\oplus y)\oplus z=(x+y-a)\oplus z=x+y+z-2a=$ $=x+(y+z-a)-a=x\oplus(y\oplus z).$ Consequently, the operation $\oplus$ is associative.

Let $x,y,z\in S_r(a)$. Then we have
\begin{eqnarray*}
(x\odot y)\odot z&=&(r(x-a)(y-a)+a)\odot z\\
&=&r^2(x-a)(y-a)(z-a)+a\\
&=&r(x-a)(y\odot z-a)+a\\
&=&x\odot(y\odot z).
\end{eqnarray*}
This means that the operation $\odot$ is associative.

Let us search for an element ${\bf 0}\in V_r(a)$ satisfying the equation ${\bf 0}\oplus x=x$ for arbitrary $x\in V_r(a)$.
So we get ${\bf 0}\oplus x={\bf 0}+x-a=x$. From this, we have an unique element ${\bf 0}=a\in V_r(a)$.

Let us search for an element ${\bf 1}\in S_r(a)$ satisfying the equation ${\bf 1}\odot x=x$ for arbitrary $x\in S_r(a)$.
So we get ${\bf 1}\odot x=r({\bf 1}-a)(x-a)+a=x$. From this, we have an unique element ${\bf 1}={1\over r}+a\in S_r(a)$.

One can check that
$$
(2a-x)\oplus x=2a-x+x-a=a={\bf 0},\ \ \ \forall x\in V_r(a),
$$
and

$$
\frac{1}{r^2(x-a)}\odot x=r(x^{-1}-a)(x-a)+a={1\over r}+a={\bf 1},\ \ \ \forall x\in S_r(a).
$$

Then keeping in mind commutativity of $\oplus$ and $\odot$
we infer that
$(V_r(a), \oplus)$ and $(S_r(a), \odot)$ are compact abelian groups.

Note that, topological groups are combination of groups and topological spaces,
i.e., they are groups and topological spaces at the same time,
such that the continuity condition for the group operations connects these two structures together and $\mathbb Q_p$ is a metric space.

Let $x_n$ and $y_n$ be two converging sequences in $V_r(a)$ and $x_n\to x$, $y_n\to y$. Then $$\lim_{n\to\infty}(x_n\oplus y_n)=\lim_{n\to\infty}(x_n+y_n-a)=\lim_{n\to\infty}x_n+\lim_{n\to\infty}y_n-a=x+y-a=x\oplus y.$$
Also, $$\lim_{n\to\infty}(-x_n)=\lim_{n\to\infty}\left(2x_n-a\right)=2x-a=-x.$$

Let $x_n$ and $y_n$ be two converging sequences in $S_r(a)$ and $x_n\to x$, $y_n\to y$. Then $$\lim_{n\to\infty}(x_n\odot y_n)=r\left(\lim_{n\to\infty}x_n-a\right)\left(\lim_{n\to\infty}y_n-a\right)+a=x\odot y.$$
Also, $$\lim_{n\to\infty}x^{-1}_n=\lim_{n\to\infty}\left({1\over{r^2(x_n-a)}}+a\right)={1\over{r^2(x-a)}}+a=x^{-1}.$$
This complete the proof.
\end{proof}

\begin{rk}\label{r1}
It should not be forgotten that the multiplication (\ref{odot}) on a sphere depends on the radius and center of this sphere.
Since the radius of a ball or sphere in $Q_p$ is a rational number in the form $r=p^m$, $m\in\Z$, the number $r$ was considered as a $p$-adic number in the determination of the product (\ref{odot}).
\end{rk}

Let $\left(V_{r_1}(a_1), \oplus_1\right)$, $\left(V_{r_2}(a_2), \oplus_2\right)$ and $\left(S_{r_1}(a_1), \odot_1\right)$, $\left(S_{r_2}(a_2), \odot_2\right)$ be the groups.

\begin{thm}\label{git} For any $r_1>0$, $r_2>0$ and $a_1, a_2\in\Q_p$ the following hold
\begin{itemize}
\item[1.] $\left(V_{r_1}(a_1), \oplus_1\right)\cong\left(V_{r_2}(a_2), \oplus_2\right)$;
\item[2.] $\left(S_{r_1}(a_1), \odot_1\right)\cong\left(S_{r_2}(a_2), \odot_2\right)$.
\end{itemize}
\end{thm}

\begin{proof}
1. Let $h:V_{r_1}(a_1)\to V_{r_2}(a_2)$ be a bijection of the following form
\begin{equation}\label{isom}
h(x)={{r_1(x-a_1)}\over{r_2}}+a_2.
\end{equation}
Then for all $x,y\in V_{r_1}(a_1)$ we get $$h(x\oplus_1 y)={{r_1(x\oplus_1 y-a_1)}\over{r_2}}+a_2={{r_1(x+y-2a_1)}\over{r_2}}+a_2={{r_1[(x-a_1)+(y-a_1)]}\over{r_2}}+a_2=$$
$$={{r_1(x-a_1)}\over{r_2}}+a_2+{{r_1(y-a_1)}\over{r_2}}+a_2-a_2=h(x)+h(y)-a_2=h(x)\oplus_2 h(y).$$
Consequently, $\left(V_{r_1}(a_1), \oplus_1\right)\cong\left(V_{r_2}(a_2), \oplus_2\right)$.

2. Let $h:S_{r_1}(a_1)\to S_{r_2}(a_2)$ be a bijection of the form (\ref{isom}).
Then for all $x,y\in S_{r_1}(a_1)$ we get $$h(x\odot_1 y)={{r_1(x\odot_1 y-a_1)}\over{r_2}}+a_2={{r_1^2(x-a_1)(y-a_1)}\over{r_2}}+a_2=$$
$$=r_2\left({{r_1(x-a_1)}\over{r_2}}+a_2-a_2\right)\left({{r_1(y-a_1)}\over{r_2}}+a_2-a_2\right)+a_2=h(x)\odot_2 h(y).$$
Consequently, $\left(S_{r_1}(a_1), \odot_1\right)\cong\left(S_{r_2}(a_2), \odot_2\right)$.
\end{proof}

\section{Haar measure on the ball and sphere}

Let $G$ be a topological group.
If $G$ is abelian and locally compact, then it is well known \cite{Kan} that it has a nonzero translation-invariant measure $\mu$, which is unique up to scalar. This is called the {\it Haar measure}.

In the field of $p$-adic numbers, one usually considers the Haar measure either on the whole space, or on its maximal ideal $\Z_p$.
That is, $X=\Q_p$ or $X=\Z_p$, and let $\Sigma$ be the minimal $\sigma$-algebra containing all open and closed (clopen) subsets of $X$.
A measure $\mu(V_{\rho})=\rho, \, V_{\rho}\in \Sigma$ is usually called a Haar measure, where $V_{\rho}$ is a ball with radius $\rho$.

However, in some cases, the problem of studying the dynamical system of a function that mapping a compact subset of $\Q_p$ to itself arises. At this time, is needed a measure defined on $\sigma$-algebra with the unit a compact set.
If this compact set has some algebraic structure, then can we look at the natural Haar measure?
If the compact set under consideration is a ball or a sphere, the answer to this question is positive, which is given in the following theorem.

Let $V_r(a)$ be the ball ($S_r(a)$ be the sphere) with the center at the point $a\in\Q_p$ and $\mathcal B$ is the algebra generated by clopen
subsets of $V_r(a)$ ($S_r(a)$). It is known that every element of $\mathcal B$ is a union of
some balls $V_{\rho}(s)\subset V_r(a)$ ($V_{\rho}(s)\subset S_r(a)$).

\begin{thm}\label{haar}
A measure $\bar\mu:\mathcal B\rightarrow p^{\mathbb Z}$ is a
Haar measure if it is defined by $\bar\mu(V_{\rho}(s))=\rho$ for all $V_{\rho}(s)\in \mathcal B$.
\end{thm}

\begin{proof}
In the previous section, it was shown that any ball and sphere is a topological abelian group with respect to the operations (\ref{oplus}) and (\ref{odot}) respectively. Then to prove this theorem, it is sufficient to prove the equalities
$$\bar\mu\left(x\oplus V_{\rho}(s)\right)=\bar\mu\left(V_{\rho}(s)\right),$$
for arbitrary $x\in V_r(a)$ and $V_{\rho}(s)\in\mathcal B$, where $x\oplus A:=\{x\oplus y| \, y\in A\}$, and
$$\bar\mu\left(x\odot V_{\rho}(s)\right)=\bar\mu\left(V_{\rho}(s)\right),$$
for arbitrary $x\in S_r(a)$ and $V_{\rho}(s)\in\mathcal B$, where $x\odot A:=\{x\odot y| \, y\in A\}$.

Let us prove the second of these two equalities, the first being easily proved by the same reasoning.

Let $x\in S_r(a)$ and $y\in V_{\rho}(s)\in\mathcal B$. So $|y-s|_p\leq\rho$, and $$|x\odot y-x\odot s|_p=|r(x-a)(y-a)+a-(r(x-a)(s-a)+a)|_p=$$ $$=|r|_p|x-a|_p|(y-a)-(s-a)|_p=|y-s|_p\leq\rho.$$
Then we conclude that $\bar\mu\left(x\odot V_{\rho}(s)\right)=\bar\mu\left(V_{\rho}(x\odot s)\right)=\bar\mu\left(V_{\rho}(s)\right)$.
\end{proof}

It is natural that the measure $\bar\mu$ introduced in the above theorem is bounded because it is introduced in a compact set, but it is not a probability measure. Therefore, in the next section we consider the normalized Haar measure
 $\mu$, defined by
 $$\mu(V_{\rho}(s))={{\bar\mu(V_{\rho}(s))}\over{\bar\mu(S_r(a))}}={{p\rho}\over{(p-1)r}}.$$

\section{Dynamical system of isometry on the invariant sphere}
It is known that given a metric space, an isometry is a transformation that maps elements to the same or another metric space, and the distance between the elements of the image in the new metric space is equal to the distance between the elements in the original metric space.

In this section, we consider the dynamical system $(S_r(a), f, \mu)$, where $f:S_r(a)\to S_r(a)$ be an isometry and $\mu$ is a normalized Haar measure.

Note that, the isometry $f:S_r(a)\to S_r(a)$ is bijection. Indeed, injectivity and continuity easily follows from isometry.
Suppose $f(S_r(a))\neq S_r(a)$. Then there exists an element $z\in S_r(a)$ such that $z\notin f(S_r(a))$. Since the function $f$ is continuous and $S_r(a)$ is compact, we have $f(S_r(a))$ is compact too. Then $S_r(a)\setminus f(S_r(a))$ is open and for this $z$, we can find an open neighbourhood $$U_{\varepsilon}(z)\subset S_r(a)\setminus f(S_r(a)).$$
Note that the field of $p$-adic numbers $\Q_p$ has a natural hierarchical structure: each ball consists of a finite number of non-intersecting balls of smaller radius\cite{Kat}. From this fact we can represented the sphere $S_r(a)$ by open balls with radii less than $\varepsilon$. Then there exists some minimal $N$, such that $$S_r(a)=\bigcup_{j=1}^N U_j,$$ where $U_j$ is open balls with radii less than $\varepsilon$. Since $z\in S_r(a)$, there exists $m\in\{1,2,...,N\}$, such that  $z\in U_m$, then $U_m\cap f(S_r(a))=\emptyset$. So $f(S_r(a))$ is covered by $N-1$ open balls, but $\bigcup^{N-1}_{j=1}f^{-1}(Uj)$ is a cover of $S_r(a)$, and all the radii are less than $\varepsilon$, because $f$ was an isometry. This contradicts the minimality of $N$. From this $f$ is surjective, also, bijective.

\begin{rk}\label{-1}
It can be seen from the above considerations that in the study of the ergodicity of the dynamics of isometric mapping on invariant spheres, the equation $T^{-1}(E)=E$ in the Definition \ref{erg} can be replaced to the equation $T(E)=E$.
\end{rk}

Let us now consider questions about the properties of isometric dynamics.

\begin{lemma}\label{mp}
Let $f:S_r(a)\to S_r(a)$ be an isometry and $\mu$ be a normalized Haar measure on $S_r(a)$. Then for any ball $V_\rho(s)\subset S_r(a)$ the following equality holds $$\mu(f(V_{\rho}(s)))=\mu(V_{\rho}(s)).$$
\end{lemma}
\begin{proof} The proof of this lemma follows easily from the fact that $f$ is an isometry. \end{proof}

\begin{pro}\label{dyn}
Let $f:S_r(a)\to S_r(a)$ be an isometry. Then for any initial point $x\in S_r(a)$ (except fixed point) the orbit
$\{f^n(x)| \, n\in \N\}$ isn't convergent.
\end{pro}

\begin{proof}
For any $x\in S_r(a)$ we have $|f(x)-x|_p>0$.
Since $f$ is isometry, we have $$|f^{n+1}(x)-f^n(x)|_p=|f(x)-x|_p=const$$ for any $n\in \N$. Therefore, the sequences $f^n(x)$ is not convergent.
\end{proof}

In \cite{LRS},\cite{RS20},\cite{RS2}-\cite{S2} the dynamics of various rational functions were considered. In these papers, we have shown that dynamics is an isometry on invariant spheres. We also proved that for any $x$ in an invariant sphere with radius $r$, the value of $|f(x)-x|_p$ does not depend on the choice of the point $x$, but depends only on the radius $r$.
In this paper we assume that for isometry $f:S_r(a)\to S_r(a)$ the value of $|f(x)-x|_p$ does not depend on the choice of the point $x$, but depends only on the radius $r$, and we denote $\rho(r)=|f(x)-x|_p$ for all $x\in S_r(a)$.

\begin{pro}\label{min}
Let $f:S_r(a)\to S_r(a)$ be an isometry and $\rho(r)=|f(x)-x|_p$ for all $x\in S_r(a)$.
Then the ball with radius $\rho(r)$ is minimal invariant ball for $f$.
\end{pro}

\begin{proof}
Since $|f(s)-s|_p=\rho(r)$, we have $f(s)\in V_{\rho(r)}(s)$ and $s\in V_{\rho(r)}(f(s))$.
By Lemma \ref{mp} and the fact that any point of a ball is its center we get
$$f(V_{\rho(r)}(s))=V_{\rho(r)}(f(s))=V_{\rho(r)}(s).$$

Let $V_{\theta}(s)\subset S_r(a)$ and  the ball $V_{\theta}(s)$
is an invariant for $f$, then $f(s)\in V_{\theta}(s)$, i.e.
$|f(s)-s|_p\leq\theta$. We have $\rho(r)=|f(s)-s|_p$ for every $s\in
S_r(a)$. So then $\rho(r)\leq\theta$. Consequently,
the ball with radius $\rho(r)$ is minimal invariant ball for $f$.
\end{proof}

\begin{thm}\label{t9}
Let $f:S_r(a)\to S_r(a)$ be an isometry and $\rho(r)=|f(x)-x|_p$ for all $x\in S_r(a)$. If the dynamical system $(S_r(a), f)$ has $k$-periodic orbit $y_0\to y_1\to ...\to y_k\to y_0$,
then the following statements hold:
\begin{itemize}
\item[1.] $y_i\in V_{\rho(r)}(y_0)$ for all $i\in\{1, 2, ..., k\}$;
\item[2.] Character of periodic points is indifferent;
\item[3.] If $\rho\leq\rho(r)$, then we have $f(S_{\rho}(y_i))\subset S_{\rho}(y_{i+1})$ for any $i\in\{0,1,...k-1\}$ and\\
$f(S_{\rho}(y_k))\subset S_{\rho}(y_0).$
\end{itemize}
\end{thm}

\begin{proof}
1. It follows from the Proposition \ref{min}.

2. Since $f$ is isometry, for all $i\in\{0,1,...,k\}$ we have the following
$$|f'(y_i)|_p=\lim_{x\to y_i}\left|{f(x)-f(y_i)}\over{x-y_i}\right|_p=1.$$ From this we conclude that
the character of periodic points $y_i$ is indifferent.

3. The proof of this assertion follows easily from the fact that $f$ is an isometry.
\end{proof}

Now we present the result about the ergodicity of the dynamical system of isometry.

\begin{thm}\label{Erg}
Let $f:S_r(a)\to S_r(a)$ be an isometry and $\rho(r)=|f(x)-x|_p$ for all $x\in S_r(a)$. Then
\begin{itemize}
\item[a)] the dynamical system $(S_r(a), f, \mu)$ is not ergodic for all $p\geq 3$;
\item[b)] the dynamical system $(S_r(a), f, \mu)$ may be ergodic if and only if $r=2\rho(r)$ for $p=2$.
\end{itemize}
Here $\mu$ is normalized Haar measure on $S_r(a)$.
\end{thm}

\begin{proof}
Note that the normalized Haar measure defined on the $\sigma$-algebra $\mathcal B$, also, every element of $\mathcal B$ is a union of
some balls $V_{\rho}(s)\subset S_r(a)$.

By Proposition \ref{min} we have $\rho(r)$ is a radius of minimal invariant ball. From the Remark \ref{-1}, the Definition \ref{erg} and the definition of normalized Haar measure we get the dynamical system $(S_r(a), f, \mu)$ may be ergodic if and only if
\begin{equation}\label{rhor}
{{p\rho(r)}\over{r(p-1)}}=1.
\end{equation}
Note that $\rho(r)<r$. Let $r=p^n\rho(r)$ for some $n\in\N$. Then from (\ref{rhor}) we have the following
\begin{equation}\label{pmn}
p^{n-1}(p-1)=1.
\end{equation}
It is not difficult to see that the equality (\ref{pmn}) holds only when $p=2$ and $n=1$.
\end{proof}

\section*{Acknowledgements}

 The author would like to thank Professors U.A. Rozikov and O.N. Khakimov for intensive discussions.

\end{document}